\documentclass[a4paper, 12pt,reqno]{amsart} 
\usepackage{amssymb}
\usepackage{amsmath}

\newtheorem{thm}{Theorem}[section]

\newtheorem{cor}{Corollary}[section]

\makeatletter
   
          \@addtoreset{equation}{section}
\makeatother

\begin{document}

\baselineskip=1.13\baselineskip   

\title[Dirac operators]{Eigenfunctions of Dirac operators
at the threshold energies}

\author{ Tomio Umeda\\}

\address{Department of Mathematical Sciences,
University of Hyogo, Himeji 671-2201, Japan
}
\email{umeda@sci.u-hyogo.ac.jp}
\date{}

\maketitle

\vspace{20pt}

\textbf{Abstract.} We show that the eigenspaces of 
the  Dirac operator $H=\alpha\cdot ( D - A(x) ) + m \beta $ at
the threshold energies $\pm m$ are
coincide with the direct sum of the zero space and the kernel of 
the Weyl-Dirac operator $\sigma\cdot ( D - A(x) )$. 
Based on this result, we describe the asymptotic limits of 
the eigenfunctions of the Dirac operator corresponding to
these threshold energies. Also, we discuss the set of vector
potentials for which the kernels of $H\mp m$ are non-trivial,
i.e. $\mbox{Ker}(H\mp m) \not = \{ 0 \}$.

 \vspace{15pt}


 \vspace{15pt}
\textbf{The 2000 Mathematical Subject Classification:} 35Q40, 35P99, 81Q10

\newpage

\section{Introduction}

This note\footnote{This note is based on joint work with 
Professor Yoshimi Sait\={o}, University of Alabama at Birmingham, USA.} 
is concerned with eigenfunctions at the threshold energies
of Dirac operators with positive mass. More precisely, the Dirac
operator which we shall deal with is of the form
\begin{equation} \label{eqn:1-1}
H=H_0+Q= \alpha \cdot D  + m \beta + Q(x), \quad D=\frac{1}{\, i \,} \nabla_x,
\,\,\, x \in {\mathbb R}^3.
\end{equation}
Here $\alpha= (\alpha_1, \, \alpha_2, \, \alpha_3)$ is
the triple of  $4 \times 4$ Dirac matrices
\begin{equation*}\label{eqn:1-2}
\alpha_j = 
\begin{pmatrix}
 \mathbf 0 &\sigma_j \\ \sigma_j &   \mathbf 0 
\end{pmatrix}  \qquad  (j = 1, \, 2, \, 3)
\end{equation*}
with  the $2\times 2$ zero matrix $\mathbf 0$ and
the triple of  $2 \times 2$ Pauli matrices
\begin{equation*}\label{eqn:1-3}
\sigma_1 =
\begin{pmatrix}
0&1 \\ 1& 0
\end{pmatrix}, \,\,\,
\sigma_2 =
\begin{pmatrix}
0& -i  \\ i&0
\end{pmatrix}, \,\,\,
\sigma_3 =
\begin{pmatrix}
1&0 \\ 0&-1
\end{pmatrix},
\end{equation*}
and
\begin{equation*}\label{eqn:1-4}
\beta=
\begin{pmatrix}
I_2& \mathbf 0 \\ \mathbf 0 & -I_2
\end{pmatrix}.
\end{equation*}
The constant $m$ is assumed to be positive.

Throughout this note we assume that  $Q(x)$ is a $4\times 4$ Hermitian
matrix-valued function. In addition to this, we shall later need
several different assumptions on
$Q(x)$ under which the operator $Q=Q(x) \times$
becomes relatively compact with respect to the operator $H_0$. 
Under these assumptions, the essential spectrum of the Dirac 
operator $H$ is given by the union of the intervals 
$(-\infty, \, -m]$ and $[m, \, +\infty)$:
\begin{equation*} \label{eqn:1-5}
\sigma_{\rm{ess}}(H) =  (-\infty, \, -m] \cup  [m, \, +\infty) .
\end{equation*}
This fact implies that
 the discrete spectrum of $H$ is contained in the 
spectral gap $(-m, \, m)$:
\begin{equation*}  \label{eqn:1-6}
\sigma_{\rm{dis}}(H)  \subset (-m , \, m).
\end{equation*}
In other words,
discrete eigenvalues with finite multiplicity
may exist 
in the spectral gap.

By the threshold energies of $H$, we mean  the values $\pm m$, the edges of the
essential spectrum $\sigma_{\rm{ess}}(H)$.
These values are normally excluded in scattering theory.
However, they are of particular  importance and of interest from the physics point of 
view. See Pickl and D\"urr \cite{PicklDurr} and Pickl \cite{Pickl}.

The aim of this note is to investigate asymptotic behaviors of
(square-integrable) eigenfunctions corresponding to the threshold
energies $\pm m$. From the mathematical point of view,  
the values $\pm m$ are critical in the following sense:
 eigenfunctions corresponding to
 eigenvalues in the spectral gap decrease rapidly at infinity;
 on the contrary, generalized eigenfunctions corresponding
 to energies in the intervals 
 $(-\infty, \, -m) \cup  (m, \, +\infty)$ 
  behave like plane waves at infinity,
 hence stay away from zero.
 For this criticality, it is interesting 
 and important, from the mathematical point
of view as well, 
 to examine the asymptotic behaviors of 
 eigenfunctions 
 corresponding to the threshold
energies $\pm m$.

A closely related question to the aim mentioned above
 is the one about 
the existence of $Q$'s which yield eigenfunctions
 of  
the operators $H$ at the threshold energies $\pm m$. 
Our answers to this question are the same as those
which were given to the question
about the existence of magnetic fields  
giving rise to zero modes
for Weyl-Dirac operators $\sigma \cdot (D -A(x))$: 
see Adam, Muratori and
Nash
\cite{AdamMuratoriNash1},
\cite{AdamMuratoriNash2}, \cite{AdamMuratoriNash3},
Balinsky and Evans \cite{BalinEvan1},
\cite{BalinEvan2},
\cite{BalinEvan3}, and
Elton \cite{Elton}.
Namely, there exist infinitely many $Q$'s which yield 
eigenfunctions
 of
the operators $H$ at the threshold energies $\pm m$, but 
the set of such $Q$'s is still rather sparse
in a certain sense.

We should note that
one can regard the operator (\ref{eqn:1-1}) as 
a generalization of 
the Dirac operator of the form
\begin{equation} \label{eqn:1-7}
\alpha\cdot \big(D - A(x) \big) + m \beta + q(x) I_4 , 
\end{equation}
where $(q, A)$ is an electromagnetic potential,
by  taking $Q(x)$ to be $- \alpha\cdot A(x) + q(x) I_4$.
To formulate the main results of the present
note, we shall need to deal with
the operator (\ref{eqn:1-7}) in the case
where $m=0$ and $q(x)\equiv 0$.
In this case, the operator
(\ref{eqn:1-7}) becomes of the form
\begin{equation*}   \label{eqn:1-7-1}
\alpha\cdot \big(D - A(x) \big)
=
\begin{pmatrix}
 \mathbf 0 &\sigma \cdot (D -A(x)) \\ 
\sigma \cdot (D -A(x)) & \mathbf 0 
\end{pmatrix}.
\end{equation*}
In this way, the Weyl-Dirac operator
\begin{equation}  \label{eqn:1-wd}
 T = \sigma \cdot (D -A(x)) 
\end{equation}
 mentioned above naturally appears in our setting.

Also, we should like to note that 
the operator (\ref{eqn:1-1}) generalizes
the Dirac operator of the form
\begin{equation} \label{eqn:1-8}
\alpha\cdot D  + m(x) \beta  + q(x) I_4, 
\end{equation}
where $m(x)$, considered to be a variable mass,
converges to a positive
constant $m_{\infty}$ at infinity in an appropriate manner.

Spectral properties of the operator (\ref{eqn:1-8})
have been extensively studied under various
assumptions on $m(x)$ and $q(x)$ in recent years.
See Kalf and Yamada \cite{KalfYamada},
Kalf, Okaji and Yamada \cite{KalfOkajiYamada},
Schmidt and Yamada \cite{SchmidtYamada},
Pladdy \cite{Pladdy} and
Yamada \cite{Yamada}.

\vspace{15pt}

\noindent
\textbf{Notation.}

\noindent
By $L^2 =L^2({\mathbb R}^3)$, we mean the Hilbert space of
square-integrable functions on ${\mathbb R}^3$, and 
we introduce  a Hilbert space ${\mathcal L}^2$ by
     ${\mathcal L}^2 = [L^2({\mathbb R}^3)]^4$, where 
the inner product  
is given by
\begin{equation*}
    (f, g)_{{\mathcal L}^2}
 = \sum_{j=1}^4 (f_j, g_j)_{L^2}
\end{equation*}
for  $f = {}^t(f_1, f_2, f_3, f_4)$ 
and
$g = {}^t(g_1, g_2, g_3, g_4)$.

By $L^{2, s}({\mathbb R}^3)$, we mean the weighted $L^2$ space
defined by
\begin{equation*}
 L^{2, s}({\mathbb R}^3)
  :=\{ \, u \; | \; \langle x \rangle^s u 
 \in L^2({\mathbb R}^3) \, \}
\end{equation*}
with the inner product
\begin{equation*}
 (u, \, v)_{L^{2, s}}
:=
\int_{{\mathbb R}^3} 
   \langle x \rangle^{2s} u(x)  \, \overline{v(x)} \, dx, 
\end{equation*}
where
\begin{equation*}  \label{eqn:bracket-x}
\langle x \rangle = \sqrt{1 + |x|^2 \,}.
\end{equation*}
We introduce  the Hilbert space
     ${\mathcal L}^{2,s} = [L^{2,s}({\mathbb R}^3)]^4$ with
the inner product  
\begin{equation*}
    (f, g)_{{\mathcal L}^{2,s}}
 = \sum_{j=1}^4 (f_j, g_j)_{L^{2,s}}.
\end{equation*}

By $H^{1}({\mathbb R}^3)$ we denote
the Sobolev space of order $1$,
and by  
${\mathcal H}^{1}$ we mean the
Hilbert space $[H^{1}({\mathbb R}^3)]^4$.
By $S({\mathbb R}^3)$, we mean
 the Schwartz class of rapidly decreasing
functions on ${\mathbb R}^3$, and we set 
$\mathcal S = [S({\mathbb R}^3)]^4$.

\vspace{15pt}
When we mention the Weyl-Dirac operator 
$\, T= \sigma \cdot (D -A(x)) \,$,
we must handle
 two-vectors (two
components spinors)
which will be denoted by $\varphi$, $\psi$, etc.
Note that the Hilbert space for the Weyl-Dirac
operator is 
$[L^2({\mathbb R}^3)]^2$.

\vspace{15pt}


\section{Massless Dirac operators}

In this section, we shall treat the Dirac operator
(\ref{eqn:1-1}) in the massless case $m=0$
under Assumption (Q) below. Namely,
we shall consider the operator
\begin{equation} \label{eqn:2-1}
H=H_0+Q= \alpha \cdot D  + Q(x).
\end{equation}
We need discussions 
about the operator (\ref{eqn:2-1}) 
in order to formulate the main results
of this note, 
which will be stated in
section  \ref{sec:massive}.

\newpage

\noindent
\textbf{Assumption (Q).} 

\noindent
Each element $q_{jk}(x)$ 
($j, \, k =1, \, \cdots, \, 4$) of $Q(x)$ is 
a measurable function satisfying
\begin{equation} \label{eqn:2-2}
| q_{jk}(x) | \le C_{\! q} \langle x \rangle^{-\rho} 
\quad   ( \rho >1 )
\end{equation}
where $C_{\! q}$ is a positive constant.

\vspace{15pt}

One should note that, under Assumption(Q), the Dirac operator
(\ref{eqn:2-1}) is a self-adjoint operator in $\mathcal L^2$ with  
$\mbox{Dom}(H) = {\mathcal H}^1$.
The self-adjoint realization will be denoted 
by $H$ again. 
With an abuse of notation, we shall write $Hf$  
\textit{in the distributional sense} for
$f \in {\mathcal S}^{\prime}$ 
whenever it makes sense.

\vspace{15pt}

\noindent
D{\scriptsize EFINITION.} \ 
By a zero mode, we mean a function 
$f \in \mbox{Dom}(H)$ which satisfies
\begin{equation*}  \label{eqn:2-3}
Hf=0.
\end{equation*}
By a zero resonance, we mean a function 
$f \in {\mathcal L}^{2, -s}\setminus \mathcal L^2$, for
some $s \in (0, \, 3/2]$,
 which satisfies
$\, Hf=0 \,$ in the distributional sense.

\vspace{15pt}

It is evident that a zero mode of $H$ is an eigenfunction of $H$
corresponding to the eigenvalue $0$, i.e.,
a zero mode is
 an element of $\mbox{Ker}(H)$,
 the kernel of the self-adjoint
operator $H$. 

\vspace{10pt}

We now  state  results
which will be needed in section \ref{sec:massive}.

\vspace{10pt}

\begin{thm} \label{thm:th-lim}
Suppose Assumption {\rm(Q)} is satisfied. Let $f$ be  
a zero mode of the operator {\rm(\ref{eqn:2-1})}.
Then for any $\omega \in {\mathbb S}^2$ 
\begin{equation}   \label{eqn:lim-1}
\lim_{r\to +\infty}r^{2}f(r\omega)
= 
-\frac{i}{\, 4\pi \,} \,
(\alpha \cdot \omega) \! \!
\int_{{\mathbb R}^3}   Q(y) f(y) \, dy,
\end{equation}
where the convergence is uniform
with respect to $\omega \in {\mathbb S}^2$.
\end{thm}

\vspace{10pt}

In connection with the expression $f(r\omega)$ in (\ref{eqn:lim-1}), 
it is worthy to
note that every zero mode is a continuous function.
This fact was shown in \cite{SaitoUmeda2}.

Theorem \ref{thm:resonance} below 
means that zero
resonances do not exist under
the restrictions on $\rho$ and $s$ imposed
in the theorem. 

\vspace{5pt}
\begin{thm} \label{thm:resonance}
Suppose Assumption {\rm(Q)} is satisfied with
$\rho > 3/2$. If $f$ belongs to
${\mathcal L}^{2, \, -s}$ for some 
$s$ with $0 < s \le \min\{3/2, \, \rho - 1 \}$ and
satisfies $Hf=0$ in the 
distributional sense, then 
$f \in {\mathcal H}^1$.
\end{thm}

\vspace{5pt}
For the proofs of Theorems \ref{thm:th-lim} and 
\ref{thm:resonance}, see \cite{SaitoUmeda1} and
\cite{SaitoUmeda2}.

\vspace{5pt}

As can be easily understood from
the discussions   in the introduction, 
the Dirac operator 
$H=\alpha \cdot D  + Q(x)$
in (\ref{eqn:2-1}) is a natural generalization of
the Weyl-Dirac operator
$\, T= \sigma \cdot (D -A(x)) \,$.
Accodingly, we obtain results on
the Weyl-Dirac operator as  corollaries to
Theorems \ref{thm:th-lim} and \ref{thm:resonance}.
To state these theorems, we have to make an assumption
on the vector potential $A(x)$,
in accordance with Assumption (Q).

\vspace{10pt}

\noindent
\textbf{Assumption (A1).} 

\noindent
Each element $A_{j}(x)$ 
($j =1, \, 2, \, 3$) of $A(x)$ is 
a real-valued measurable function satisfying
\begin{equation} \label{eqn:2-1WD}
| A_{j}(x) | \le C_{\! a} \langle x \rangle^{-\rho} 
\quad   ( \rho >1 )
\end{equation}
where $C_{\! a}$ is a positive constant.

\vspace{10pt}

Assumption (A1) assures that 
$T =\sigma\cdot(D-A(x))$ is a self-adjoint operator
in $[L^2({\mathbb R}^3)]^2$  
with domain $[H^1({\mathbb R}^3)]^2$.  

\vspace{10pt}

\begin{thm} \label{thm:th-limWD}
Suppose 
Assumption {\rm(A1)} is satisfied.
Let $\psi$ be  
a zero mode of the Weyl-Dirac
operator $T =\sigma\cdot(D-A(x))$.
Then for any $\omega \in {\mathbb S}^2$ 
\begin{gather}
\begin{split}   \label{eqn:lim-1WD}
{}&\lim_{r\to +\infty}r^{2}\psi(r\omega)     \\
&= 
\frac{i}{\, 4\pi \,}
\int_{{\mathbb R}^3}  
\big\{ 
  \big( \omega\cdot A(y) \big)  I_2
    + i\sigma \cdot 
    \big( \omega\times A(y) \big) 
   \big\}  \psi(y) \, dy,
\end{split}
\end{gather}
where the convergence is uniform
with respect to $\omega \in {\mathbb S}^2$.
\end{thm}

\vspace{5pt}
\begin{thm} \label{thm:WDresonance}
Suppose Assumption {\rm(A1)} is satisfied with
$\rho > 3/2$. If $\psi$ belongs to
$[{L}^{2, \, -s}({\mathbb R}^3)]^2$ for some 
$s$ with $0 < s \le \min\{3/2, \, \rho - 1 \}$ and
satisfies $T\psi=0$ in the 
distributional sense, then 
$\psi  \in [{H}^1({\mathbb R}^3)]^2$.
\end{thm}


\vspace{10pt}

\section{Dirac operators with positive mass} \label{sec:massive}

In this section, we shall restrict ourselves to the Dirac 
operators with  a vector potential
\begin{equation} \label{eqn:3-do}
H=\alpha\cdot \big(D - A(x) \big) + m \beta, 
\end{equation}
where $m>0$.

One of our main results characterizes eigenfunctions 
of the Dirac operator $H$ in (\ref{eqn:3-do}) 
at the threshold eigenvalues
$\pm m$ in terms of zero modes
of the Weyl-Dirac operator $T$ in (\ref{eqn:1-wd}).

\vspace{15pt}

\begin{thm} \label{thm:th-DO}
Suppose 
Assumption {\rm(A1)} is satisfied. Then
\vspace{3pt}
\begin{itemize}
\item[(i)] $f\in \mbox{\rm Ker}(H-m) \Longleftrightarrow
\exists \psi \in \mbox{\rm Ker}(T)$ such that
$\displaystyle{f= \binom{\psi}{0}  }$\mbox{\rm;}
\vspace{6pt}
\item[(ii)] $f\in \mbox{\rm Ker}(H+m) \Longleftrightarrow
\exists \psi \in \mbox{\rm Ker}(T)$ such that
$\displaystyle{f= \binom{0}{\psi}  }$.
\end{itemize}
\end{thm}

\vspace{15pt}

It is of some interest to point out  
that Theorem \ref{thm:th-DO} implies that eigenfunctions 
and eigenspaces corresponding to the 
threshold eigenvalues $\pm m$ are independent of $m$.
Also we point out that
Theorem \ref{thm:th-DO} implies
\begin{gather}  
\mbox{\rm Ker}(H-m)= \mbox{\rm Ker}(T)\oplus\{ 0 \},  
\label{eqn:3-ker-} \\
\noalign{\vskip 6pt}
\mbox{\rm Ker}(H+m)= \{ 0 \} \oplus \mbox{\rm Ker}(T). 
\label{eqn:3-ker+} 
\end{gather}

\vspace{8pt}

It is immediate that Theorem \ref{thm:th-DO}, together with
Theorem \ref{thm:th-limWD}, yields the following corollary.

\vspace{10pt}

\begin{cor}
Suppose 
Assumption {\rm(A1)} is verified. 
Let $u_{\psi}(\omega)$ be the continuous
function on ${\mathbb S}^2$ defined by (\ref{eqn:lim-1WD}).
\vspace{3pt}
\begin{itemize}
\item[(i)] If $f\in \mbox{\rm Ker}(H-m)$,
then $f$ is continuous on ${\mathbb R}^3$
and
satisfies that
for any $\omega \in {\mathbb S}^2$ 
\begin{equation*}   
\lim_{r\to +\infty}r^{2}f(r\omega)=\binom{u_{\psi}(\omega)}{0}, 
\end{equation*}
where the convergence is uniform
with respect to $\omega \in {\mathbb S}^2$.
\vspace{6pt}
\item[(ii)] If $f\in \mbox{\rm Ker}(H+m)$,
then
$f$ is continuous on ${\mathbb R}^3$
and
satisfies that
for any $\omega \in {\mathbb S}^2$ 
\begin{equation*}   
\lim_{r\to +\infty}r^{2}f(r\omega)=\binom{0}{u_{\psi}(\omega)}, 
\end{equation*}
where the convergence is uniform
with respect to $\omega \in {\mathbb S}^2$.
\end{itemize}
\end{cor}

The conclusions of Theorem \ref{thm:th-DO} are
valid under weaker assumptions than 
Assumption (A1). Indeed, we shall introduce 
two  assumptions, which are
quite different from each other.
 The one (Assumption (A2) below)
 needs continuity
of the vector potentials, but it allows a slightly
slower decay
at infinity.  
The other one (Assumption (A3) below) allows the vector
potentials with local singularities. 


\newpage

\noindent
\textbf{Assumption (A2).} 

\noindent
Each element $A_{j}(x)$  is 
a real-valued continuous function satisfying
\begin{equation} \label{eqn:3-2}
 A_{j}(x) = o(|x|^{-1})  \quad  (|x| \to + \infty). 
\end{equation}

\vspace{15pt}

\noindent
\textbf{Assumption (A3).} 

\noindent
Each element $A_{j}(x)$  is 
a real-valued measurable function satisfying
\begin{equation}  \label{eqn:3-3}
A_j \in L^3( {\mathbb R}^3 ).
\end{equation}

\vspace{8pt}

It is obvious that under Assumption (A2), 
$-\alpha \cdot A$ is a bounded self-adjoint 
operator in $\mathcal L^2$, hence
the Dirac operator
(\ref{eqn:3-do}) is a self-adjoint operator in $\mathcal L^2$ with  
$\mbox{Dom}(H) = {\mathcal H}^1$.

As for the Dirac operator
(\ref{eqn:3-do}) under Assumption (A3),
one can show that
$\alpha\cdot \big(D - A(x) \big) + m \beta$
is a relatively compact perturbation of
the operator $\alpha\cdot D + m \beta$.
Therefore we find that the formal expression
(\ref{eqn:3-do}) admits the self-adjoint 
realization in $\mathcal
L^2$ with  
$\mbox{Dom}(H) = {\mathcal H}^1$.

\vspace{10pt}

\begin{thm}    \label{thm:th-DO-23}
Suppose either of 
Assumption {\rm(A2)}  or  Assumption {\rm(A3)} is satisfied.
 Then all the conclusions of Theorem \ref{thm:th-DO} hold.
 \end{thm}
 
 \vspace{10pt}

We should like to mention relevant works on
the Weyl-Dirac operator $T$ by
 Balinsky and Evans, and Elton. 
In their work \cite{BalinEvan2},
Balinsky and Evans showed sparseness of the set of 
vector potentials and derived
an estimate of the dimension of 
the subspace consisting of zero modes of $T$ under 
Assumption (A3). 
In a similar spirit to \cite{BalinEvan2},
Elton \cite{Elton}
investigated, under Assumption (A2),
the local structure of 
the set of vector potentials which produce
zero modes with multiplicity $k\ge0$.

Combining the results in \cite{Elton}
with Theorem \ref{thm:th-DO-23} above,
we get Theorem \ref{thm:th-WD-El}  below.
To formulate the theorem, we need to introduce a Banach space
of  vector potentials 
\begin{equation*}   \label{eqn:A-1}
{\mathcal A} =
\{ \; A 
\; | \, A_j(x) \in C^0({\mathbb R}^3, \mathbb R),  
\; A_j(x) = o(|x|^{-1}) \;\; \mbox{as }|x| \to +\infty
\; \}
\end{equation*}
equipped with the norm
\begin{equation*}  \label{eqn:A-2}
\Vert A \Vert =  
\sup_{x \in {\mathbb R}^3} \langle x \rangle |A(x)|.
\end{equation*}

\vspace{10pt}

\begin{thm}    \label{thm:th-WD-El}
Suppose Assumption {\rm(A2)} is satisfied. Let
\begin{equation*}   \label{eqn:Z-k}
{\mathcal Z}_k^{\pm} = 
\{ \; A \in {\mathcal A} \; | \;
 \mbox{\rm dim (Ker}(H\mp m)) = k\, \}
\end{equation*}
for $k=0, \, 1, \, 2, \, \cdots$. Then
\begin{itemize}
\item[(i)] ${\mathcal Z}_k^+ ={\mathcal Z}_k^-$ for any $k$ \mbox{\rm ;}
\item[(ii)]  ${\mathcal Z}_0^{\pm}$ is an open dense subset
of $\mathcal A$
\mbox{\rm ;}
\item[(iii)] for any $k$ and any open subset 
$\Omega \subset {\mathbb R}^3$, $\Omega\not= \varnothing$,
\begin{equation*}  \label{eqn:Z-nonempty}
[C_0^{\infty}({\Omega})]^3 
 \cap 
{\mathcal Z}_k^{\pm}  \not= \varnothing.
\end{equation*}
\end{itemize}
\end{thm}

\vspace{10pt}

In a similar fashion, we combine the results in \cite{BalinEvan2}
with Theorem \ref{thm:th-DO-23} above,
and we get the following

\vspace{10pt}

\begin{thm}    \label{thm:th-WD-BE}
Suppose Assumption {\rm(A3)} is satisfied. Then
\begin{itemize}
\item[(i)] the sets 
$\big\{ 
\, A  \in [L^3({\mathbb R}^3)]^3 \; | \; 
 \mbox{\rm Ker}(H\mp m) =\{\, 0\,\} \,
\big\}$
contain
 open dense subsets of $[L^3({\mathbb R}^3)]^3$;
\item[(ii)] $\mbox{\rm dim (Ker}(H - m))=
\mbox{\rm dim (Ker}(H + m)) \le c_0 
\displaystyle{\int_{{\mathbb R}^3} |A(x)|^3 \, dx}$
for some constant $c_0$.

\end{itemize}
 \end{thm}
 
 \vspace{10pt}

The proofs of the theorems in this section will be
found in our coming paper \cite{SaitoUmeda3}.

\vspace{10pt}


\end{document}